\documentclass[10pt]{article}
\usepackage{amssymb,amsmath,amsthm}
\usepackage[numbers,sort&compress]{natbib}
\usepackage{indentfirst}
\usepackage{exscale}
\usepackage{relsize}
\usepackage{subfigure}
\usepackage{graphicx}

\newcommand{\R}{\mathbb{R}}
\theoremstyle{definition}

\theoremstyle{remark}

\numberwithin{equation}{section}
\allowdisplaybreaks
 \UseRawInputEncoding

\begin{document}
\title{\Large\bf{ Infinitely many solutions for three quasilinear Laplacian systems on weighted graphs}}
\date{}
\author {Yan Pang$^{1}$, Junping Xie$^{2}$ \footnote{Corresponding author, E-mail address: hnxiejunping@163.com}, \ Xingyong Zhang$^{1,3}$\\
{\footnotesize $^1$Faculty of Science, Kunming University of Science and Technology,}\\
 {\footnotesize Kunming, Yunnan, 650500, P.R. China.}\\
  {\footnotesize $^2$Faculty of Transportation Engineering, Kunming University of Science and Technology, Kunming, Yunnan, 650500, P.R. China.}\\
{\footnotesize $^{3}$Research Center for Mathematics and Interdisciplinary Sciences, Kunming University of Science and Technology,}\\
 {\footnotesize Kunming, Yunnan, 650500, P.R. China.}\\}

 \date{}
 \maketitle

 \begin{center}
 \begin{minipage}{15cm}
 \small  {\bf Abstract:} We investigate a generalized poly-Laplacian system with a parameter on weighted finite graph, a generalized poly-Laplacian system with a parameter and  Dirichlet boundary value on weighted locally finite  graphs, and  a $(p,q)$-Laplacian system with a parameter on  weighted locally finite  graphs. We  utilize a critical points theorem built by  Bonanno and Bisci [Bonanno, Bisci, and Regan, Math. Comput. Model. 2010, 52(1-2): 152-160], which is an abstract critical points theorem without compactness condition, to obtain that these three systems have infinitely many nontrivial solutions with unbounded norm when the parameters locate some well-determined range.
 \par
 {\bf Keywords:} infinitely many solutions, generalized ploy-Laplacian system,  $(p,q)$-Laplacian system, finite graph, locally finite graph.
\par
 {\bf 2020 Mathematics Subject Classification.} 34B15; 34B18
 \end{minipage}
 \end{center}
  \allowdisplaybreaks
 \vskip2mm
 {\section{Introduction }}
\setcounter{equation}{0}
Assume that $G=(V, E)$ is a graph, where $V$ is the vertex set and $E$ is the edge set.  $G$ is usually known as a finite graph when $V$ and $E$ have finite elements,  and $G$ is usually  known as  a locally finite graph when for any $x\in V$, there exist finite $y\in V$ satisfying $xy\in E$, where $xy$ represents an edge linking $x$ and $y$.  The weight on any given edge $xy\in E$ is denoted by $\omega_{xy}$ which is supposed to satisfy $\omega_{xy}>0$ and $\omega_{xy}=\omega_{yx}$. Moreover, we set $deg(x)=\sum_{y\thicksim x}\omega_{xy}$ for any fixed $x\in V$. Here, we use $y\thicksim x$ to represent those $y$ linked to $x$.   $d(x,y)$ represents the distance between any two points $x,y\in V$, which is defined by the minimal number of edges linking $x$ to $y$. Suppose that $\Omega$ is a subset in $V$. If there exists a positive constant $D$ such that $d(x,y)\le D$ for all $x,y\in \Omega$, then $\Omega$ is known as a bounded domain in $V$. Set
$$ \partial \Omega:=\{y \in V, y \notin \Omega: \exists x \in \Omega \text { satisfying } x y \in E\}
$$
which is known as the boundary of $ \Omega$. The interior of $\Omega$ is represented by $\Omega^{\circ}=\Omega\backslash \partial \Omega$, which obviously satisfies $\Omega^{\circ}=\Omega$.
\par
 Thereinafter, $\mu:V\rightarrow \R^+$ is supposed to be a finite measure. Set
\begin{eqnarray}
\label{Eq5}
D_{w,y}u(x):=\frac{1}{\sqrt2}(u(x)-u(y))\sqrt\frac{w_{xy}}{\mu(x)}
\end{eqnarray}
which is the directional derivative of $u:V\rightarrow\mathbb{R}$,
 and then  the gradient of $u$ is defined as
\begin{eqnarray}
\label{Eq6}
\nabla u(x):=(D_{w,y}u(x))_{y\in V}
\end{eqnarray}
that is a vector and is indexed by $y$.
Set
\begin{eqnarray}
\label{Eq7}
\Gamma(u,v)(x)=\frac{1}{2\mu(x)}\sum\limits_{y\thicksim x}w_{xy}(u(y)-u(x))(v(y)-v(x)).
\end{eqnarray}
Then it is obvious that
\begin{eqnarray}
\label{Eq8}
\Gamma(u,v)=\nabla u\cdot \nabla v.
\end{eqnarray}
Define
\begin{eqnarray}
\label{Eq10}
|\nabla u|(x)=\sqrt{\Gamma(u,u)(x)}=\left(\frac{1}{2\mu(x)}\sum\limits_{y\thicksim x}w_{xy}(u(y)-u(x))^2\right)^{\frac{1}{2}},
\end{eqnarray}
which represent  the length of $\nabla u$. Furthermore,
The length of $m$-order gradient of $u$ is represented by $|\nabla^m u|$ that is defined by
\begin{eqnarray}
\label{Eq11}
|\nabla^mu|=
 \begin{cases}
  |\nabla\Delta^{\frac{m-1}{2}}u|,& \text {if $m$ is an odd number,}\\
  |\Delta^{\frac{m}{2}}u|,&  \text {if $m$ is an even number.}
   \end{cases}
\end{eqnarray}
Here, we define $\nabla\Delta^{\frac{m-1}{2}}u$ by $(\ref{Eq6})$ with substituting $\Delta^{\frac{m-1}{2}}u$ for $u$,  and $\Delta^{\frac{m}{2}}u=\Delta(\Delta^{\frac{m}{2}-1}u)$, where the Laplacian operator $\Delta$ of $u$ is defined as
\begin{eqnarray}
\label{Eq4}
\Delta u(x):=\frac{1}{\mu(x)}\sum\limits_{y\thicksim x}w_{xy}(u(y)-u(x)).
\end{eqnarray}
\par
For any given $l>1$, set
\begin{eqnarray}\label{Eq13}
\Delta_l u(x):=\frac{1}{2\mu(x)}\sum\limits_{y\sim x}\left(|\nabla u|^{l-2}(y)+|\nabla u|^{l-2}(x)\right)\omega_{xy}(u(y)-u(x)),
\end{eqnarray}
which is known as the $l$-Laplacian operator of $u$.
 $l$-Laplacian operator obviously becomes the Laplacian operator of $u$ as $l=1$.
\par
For convenience, we set
\begin{eqnarray}\label{Eq12}
\int\limits_V u(x) d\mu=\sum\limits_{x\in V}\mu(x)u(x).
\end{eqnarray}
For any $r\in \R$ with $r\ge 1$,  set
$$
L^r(V)=\left\{u:V\to\R\Big|\int_V|u(x)|^rd\mu<\infty\right\}
$$
equipped by the norm
\begin{eqnarray}
\label{Eq2}
\|u\|_{L^r(V)}=\left(\int_V|u(x)|^rd\mu\right)^\frac{1}{r}.
\end{eqnarray}
For any $u:V\to\R$, according to the distributional sense, we write $\Delta_l$ as
\begin{eqnarray}\label{Eq14}
\int\limits_V(\Delta_l u)v d\mu=-\int\limits_V|\nabla u|^{l-2}\Gamma(u,v)d\mu,
\end{eqnarray}
where $v\in\mathcal{C}_c(V)$  and $\mathcal{C}_c(V)$ is the set of all real functions with compact support.
Furthermore,  a more general operator $\pounds_{m,l}$ could be defined as
\begin{eqnarray}\label{eq9}
\int\limits_V(\pounds_{m,l}u)\phi d\mu=
 \begin{cases}
  \int_V|\nabla^m u|^{l-2}\Gamma(\Delta^{\frac{m-1}{2}}u,\Delta^{\frac{m-1}{2}}\phi)d\mu,& \text { when $m$ is an odd number},\\
  \int_V|\nabla^m u|^{l-2}\Delta^{\frac{m}{2}}u\Delta^{\frac{m}{2}}\phi d\mu,&  \text { when $m$ is an even number}.
    \end{cases}
\end{eqnarray}
 for any $\phi\in\mathcal{C}_c(V)$,  where $l\in \R$ with $l>1$ and $m\in \mathbb N$. $\pounds_{m,p}$ is known as the poly-Laplacian of $u$  as $m=2$, and $\pounds_{m,p}$ degenerates to the $l$-Laplacian operator as $m=1$. Those above concepts and more related details refer to  \cite{Chung2005} and \cite{Yamabe 2016}.

\par
 In this paper, we focus on  the existence of infinitely many solutions for the following  generalized poly-Laplacian system on finite graph:
\begin{eqnarray}
\label{eq2}
 \begin{cases}
  \pounds_{m_1,p}u+h_1(x)|u|^{p-2}u=\lambda F_u(x,u,v),\;\;\;\;\hfill x\in V,\\
  \pounds_{m_2,q}v+h_2(x)|v|^{q-2}v=\lambda F_v(x,u,v),\;\;\;\;\hfill x\in V,\\
   \end{cases}
\end{eqnarray}
where $G=(V,E)$ is a finite graph, $m_i\in \mathbb N$, $h_i:V\rightarrow\mathbb{R}^+,i=1,2$, $1<p,q<+\infty$, $\lambda>0$, and $F:V\times\mathbb{R}^2\rightarrow\mathbb{R}$.
\par
Moreover, if $G=(V,E)$ is a locally finite graph, we focus on the existence of infinitely many solutions for the following generalized poly-Laplacian system with Dirichlet boundary condition:
\begin{eqnarray}
\label{eq3}
 \begin{cases}
  \pounds_{m_1,p}u=\lambda F_u(x,u,v),\;\;\;\;\hfill x\in \Omega^\circ,\\
  \pounds_{m_2,q}v=\lambda F_v(x,u,v),\;\;\;\;\hfill x\in \Omega^\circ,\\
  |\nabla^j u|=0,\;\;\;\;\hfill x\in \partial\Omega,\ 0\leq j\leq m_1-1,\\
  |\nabla^i v|=0,\;\;\;\;\hfill x\in \partial\Omega,\ 0\leq i\leq m_2-1,\\
   \end{cases}
\end{eqnarray}
where $1<p,q<+\infty$, $\lambda>0$, $m_i\in \mathbb N$ with $m_i\geq1$, $i=1,2$, and $\Omega\subset G(V,E)$ is a bounded domain.
\par
 Finally, we are also concerned with the existence of infinitely many solutions for the following $(p,q)$-Laplacian system on locally finite graph  $G=(V,E)$:
\begin{eqnarray}
\label{eq4}
 \begin{cases}
   -\Delta_p u+h_1(x)|u|^{p-2}u=\lambda F_u(x,u,v),\;\;\;\;\hfill x\in V,\\
   -\Delta_q v+h_2(x)|v|^{q-2}v=\lambda F_v(x,u,v),\;\;\;\;\hfill x\in V,\\
 \end{cases}
\end{eqnarray}
where $-\Delta_p$ and $-\Delta_q$ are defined by (\ref{Eq13}) with $l=p,q$, $p\ge2$ and $q\ge2$, $F:V\times \R^2 \to \R$, $h_i:V\rightarrow\mathbb{R}^+,i=1,2$,  and $\lambda> 0$.

\par
 With the development of machine learning, data analysis, social network, image processing and traffic network, the analysis on graphs has attracted some attentions (\cite{Alkama S2014,Arnaboldi V 2015,Bini A A 2015,Elmoataz A2015,Ta V T 2010,Ta V T 2008}). In particular, recently, in \cite{Yamabe 2016} and \cite{Yamabe 2017},  Grigor¡¯yan-Lin-Yang studied several nonlinear elliptic equations on graphs and first established the Sobolev spaces and the variational framework on graphs. Subsequently, there have been some works on $p$-Laplacian equations and more general poly-Laplacian equations on graphs.
For example, in \cite{Pinamonti 2022}, Pinamonti  and Stefani studied some semi-linear equations with the poly-Laplacian operator on locally finite graphs. They established some existence results of weak solutions     via a variational method  by using the continuity properties of the energy functionals.
In \cite{Shao 2023}, Shao studied a nonlinear $p$-Laplacian equation on a locally finite graph. Some existence results of positive
solutions and positive ground state solutions are  established by exploiting  the mountain pass theorem and the Nehari manifold. More related results also refer to, for example, \cite{Ge H2018}, \cite{Jiang W2018}, \cite{Imbesi 2023}, \cite{Shao 2023(2)} and \cite{Shao M 2023}.

\par
In addition to the case of single equations, recently, the study of systems on graphs has also yielded some results.
For example, in \cite{Zhang 2022}, Zhang-Zhang-Xie-Yu considered the system (\ref{eq2}) with $\lambda=1$. They supposed that $F$ takes on  the super-$(p,q)$ growth  and  then established the existence result of a nontrivial solution by exploiting the mountain pass theorem. They also established a multiplicity result by utilizing the symmetric mountain pass theorem.
In \cite{Yu 2023},  Yu-Zhang-Xie-Zhang considered (\ref{eq3}) and  system (\ref{eq4}) with $p=q$, $\lambda=1$ and $F(x,u)=-K(x,u)+W(x,u)$ for all $x\in V$. By utilizing the mountain pass theorem, they achieved that (\ref{eq3}) has a nontrivial solution.
In \cite{Yang 2023}, Yang-Zhang investigated (\ref{eq4}) with perturbations and two parameters $\lambda_1$ and $\lambda_2$. Under the assumptions that the nonlinearity satisfies a sub-$(p, q)$ conditions, they achieved that system has at least one
nontrivial solution by Ekeland¡¯s variational principle. When the nonlinearity equipped the super-$(p, q)$ conditions, they established that system has at least one nontrivial solution with positive energy and one nontrivial solution with negative energy by exploiting mountain pass theorem and Ekeland¡¯s variational principle.
In \cite{Shao 2023(2)}, when $h_1(x)=\lambda a +1$ and $h_2(x)=\lambda b +1$, Shao studied (\ref{eq4}) with $p=q$. By  Nehari manifold method and some analytical techniques, under some suitable assumptions on the potentials and nonlinear terms, they proved that  system possesses a ground state solution $(u_\lambda, v_\lambda)$ when the parameter $\lambda$ is sufficiently large.
\par
Our investigation are mainly motivated by the above mentioned works and \cite{Bonanno G2009}-\cite{Bonanno G2010(2)}. In \cite{Bonanno G2009}, Bonanno and Bisci established the existence result of a sequence $\{u_n\}$ of critical points for the functional $f_\lambda:=\Phi-\lambda\Psi$ with  $\lambda\in\mathbb{R}$, and get a
well-determined interval of the parameter $\lambda$. In \cite{Bonanno G2010(2)}, Bonanno and Bisci obtained that a class of quasilinear elliptic system in the Euclidean framework possesses infinitely many weak solutions by the abstract theorem established in \cite{Bonanno G2009}. In the present paper, we shall also apply the critical points theorem developed by  Bonanno and Bisci \cite{Bonanno G2009} to system (\ref{eq2}), (\ref{eq3}) and (\ref{eq4}), and obtain that these systems have infinitely many nontrivial solutions with unbounded norm when the parameters $\lambda$ locate some well-determined ranges. To the best of our knowledge, there seemed to be no works to investigate the existence of infinitely
many solutions for equations or systems on finite graph or locally finite graph. Our works is a preliminary attempt in this field.

\vskip2mm
{\section{Preliminaries}}
\setcounter{equation}{0}
\par
In this section, we recall some basic knowledge on the Sobolev space on graph. More details refer to \cite{Yamabe 2016,Yang 2023, Zhang 2022}. We also recall
an abstract critical point theorem built in  \cite{Bonanno G2009}, which is exploited to prove our main results.
\par
Suppose that $G=(V,E)$ is a finite graph.
For any fixed $m\in \mathbb N$ and any fixed $l\in \R$ with $l>1$, set
\begin{eqnarray*}
W^{m,l}(V)=\left\{u:V\to\R\right\}
\end{eqnarray*}
 equipped with the norm
\begin{eqnarray}
\label{Eq1}
\|u\|_{W^{m,l}(V)}=\left(\int_V(|\nabla^m u(x)|^l+h(x)|u(x)|^l)d\mu\right)^\frac{1}{l},
\end{eqnarray}
where $h(x)>0$ for all $x\in V$. $W^{m,l}(V)$ is a Banach space with finite dimension.
\par
Suppose that $G=(V,E)$ is a locally finite graph and $\Omega$ is a bounded domain in $V$. For any fixed $l\in \R$ with $l>1$ and any fixed $m\in\mathbb N$, set
\begin{eqnarray*}
W^{m,l}(\Omega)=\left\{u:\Omega\to\R\right\}
\end{eqnarray*}
equipped with the norm
\begin{eqnarray*}
\label{Eq2}
\|u\|_{W^{m,l}(\Omega)}=\left(\sum_{k=0}^m\int_{\Omega\cup\partial\Omega}(|\nabla^k u(x)|^ld\mu\right)^\frac{1}{l}.
\end{eqnarray*}
\par
Define
\begin{eqnarray*}
C_0^m(\Omega)=\{u:\Omega\rightarrow\mathbb{R}|u=|\nabla u|=\cdots=|\nabla^{m-1} u|=0\}.
\end{eqnarray*}
$W_0^{m,l}(\Omega)$ is seen as the completion of $C_0^m(\Omega)$ in $W^{m,l}(\Omega)$. $W_0^{m,l}(\Omega)$ is a finite dimensional Banach space since $\Omega$ is a finite set. On $W_0^{m,l}(\Omega)$, one can also equip the following norm
\begin{eqnarray*}
\label{Eq2}
\|u\|_{W_0^{m,l}(\Omega)}=\left(\int_{\Omega\cup\partial\Omega}(|\nabla^k u(x)|^ld\mu\right)^\frac{1}{l},
\end{eqnarray*}
Then $\|u\|_{W_0^{m,l}(\Omega)} $ is equivalent to $\|u\|_{W^{m,l}(\Omega)}$.
\par
Suppose that $G=(V,E)$ is a locally finite graph. $W^{1,l}(V)$ $(l>1)$ is the completion of $\mathcal{C}_c(V)$ based on the norm
\begin{eqnarray*}
\|u\|_{W_h^{1,l}(V)}=\left(\int_V(|\nabla u(x)|^l+h(x)|u(x)|^l)d\mu\right)^\frac{1}{l},
\end{eqnarray*}
where $h:V\to \R$ and there exists a positive constant $h_0$ such that $h(x)\geq h_0 $. Set the space
\begin{eqnarray*}
W_h^{1,l}(V)=\left\{u\in W^{1,l}(V)\big|\int_V h(x)|u(x)|^ld\mu <\infty\right\}
\end{eqnarray*}
equipped with the norm $\|u\|_{W_h^{1,l}(V)}$.

  \vskip2mm
\noindent
{\bf Lemma 2.1.} (\cite{Yamabe 2016, Zhang 2022}) {\it Suppose that $G=(V,E)$ is a finite graph. For any $\psi\in W^{m,l}(V)$, there exists
$$\|\psi\|_{\infty,V}\leq K_l\|\psi\|_{W^{m,l}(V)},$$
where $\|\psi\|_{\infty}=\max_{x\in V}|\psi(x)|$ and $K_l=\left(\frac{1}{\mu_{\min}h_{\min}}\right)^\frac{1}{l}$ with $\mu_{\min}=\min_{x\in V} \mu(x)$ and $h_{\min}=\min_{x\in V} h(x)$.}

 \vskip2mm
\noindent
{\bf Lemma 2.2.} (\cite{Yamabe 2016,Zhang 2022}) {\it Suppose that $G=(V,E)$ is a locally finite graph  and $\Omega$ be a bounded domain in $V$ satisfying $\Omega^\circ\neq\emptyset$. Let $m\in \mathbb N$ and $l>1$. Then $W_0^{m,l}(\Omega)$ is continuously embedded into $L^\theta(\Omega)$ for all $1\leq \theta\leq +\infty$. In particular, there exists a positive constant $C(m,l,\Omega)$ which just depends on $m,l$ and $\Omega$ satisfying
$$\left(\int_\Omega|u|^qd\mu\right)^\frac{1}{q} \leq C(m,l,\Omega) \left(\int_{\Omega\cup\partial\Omega}|\nabla^mu|^pd\mu\right),$$
$$\|u\|_{\infty,\Omega} \leq \frac{C}{\mu_{\min,\Omega}^{1/l}}\|u\|_{W_0^{m,l}(\Omega)},$$
for all $1\leq \theta\leq +\infty$ and all $u\in W_0^{m,l}(\Omega)$, where $\|u\|_{\Omega,\infty}=\max_{x\in\Omega}|u(x)|$ and $\mu_{\min,\Omega}=\min_{x\in \Omega}\mu(x)$. Moreover, $W_0^{m,l}(\Omega)$ is pre-compact, that is, if $\{u_n\}$ is bounded in $W_0^{m,l}(\Omega)$, then up to a subsequence, there exists some $u\in W_0^{m,l}(\Omega)$ such that $u_n\rightarrow u$ in $W_0^{m,l}(\Omega)$.
}
 \vskip2mm
\noindent
{\bf Lemma 2.3.} (\cite{Yang 2023}) {\it Suppose that $G=(V,E)$ is a locally finite graph, and  $h(x)>h_0$ and $\mu(x)>\mu_0$ for all $x\in V$, some $h_0>0$ and some $\mu_0>0$. If $(H_1)$ holds, then $W_h^{1,l}(V)$ is continuously embedded into $L^r(V)$ for all $1<l\leq r\leq \infty$, and the following inequalities hold:
\begin{eqnarray*}\label{Eq16}
\|u\|_\infty\leq \frac{1}{h_0^{1/l}\mu_0^{1/l}}\|u\|_{W_h^{1,l}(V)}
\end{eqnarray*}
and
\begin{eqnarray*}\label{Eq17}
\|u\|_{L^r(V)}\leq \mu_0^\frac{l-r}{lr}h_0^{-\frac{1}{l}}\|u\|_{W_h^{1,l}(V)}\;\;\mbox{for all}\; l\leq r<\infty.
\end{eqnarray*}}

 \noindent
{\bf Lemma 2.4.} (\cite{Bonanno G2009}) {\it Assume that $X$ is a reflexive real Banach space, $\Phi,\Psi:X\to\mathbb{R}$ are two G$\hat{a}$teaux differentiable functional satisfying  $\Phi$ is continuous, sequentially weakly lower semicontinuous and  coercive,  and $\Psi$ is sequentially weakly upper semicontinuous. For each $r>\inf_X\Phi$, set
 $$\varphi(r):=\inf_{u\in\Phi^{-1}([-\infty,r])}\frac{\left(\sup_{u\in\Phi^{-1}([-\infty,r])}\Psi(u)\right)-\Psi(u)}{r-\Phi(u)}$$
 and
 $$\gamma:=\liminf_{r\rightarrow+\infty}\varphi(r).$$
 Then,\\
 $(a)$  if $\gamma<+\infty$, for each $\lambda\in(0,\frac{1}{\gamma})$, the following alternative holds: either\\
   \indent $(a_1)$ $I_\lambda:=\Phi-\lambda\Psi$ admits a global minimum, or\\
   \indent $(a_2)$ there exists a sequence $\{u_n\}$ of critical points (local minima) of $I_\lambda$ satisfying $\lim_{n\rightarrow\infty}\Phi(u_n)=+\infty.$\\
 $(b)$ if $\delta<+\infty$, for each $\lambda\in(0,\frac{1}{\delta})$, the following alternative holds: either\\
   \indent $(b_1)$  there exists a global minimum of $\Phi$ that is a local minimum of $I_\lambda$, or\\
   \indent $(b_2)$  there exists a sequence of pairwise distinct critical points (local minima) of $I_\lambda$ that weakly converges to a global minimum of $\Phi$. }

\vskip2mm
{\section{Result and proofs for  system (\ref{eq2})}}
  \setcounter{equation}{0}
In this section, we investigate the generalized poly-Laplacian system (\ref{eq2}) and obtain the following result.
\par
Let
\begin{eqnarray}
\label{3.1}
\varrho_V=\max\left\{\frac{1}{p}\int_Vh_1(x)d\mu,\frac{1}{q}\int_Vh_2(x)d\mu\right\}, \ \ K_V=\max\left\{\frac{1}{\mu_{\min}h_{\min}},\frac{1}{\mu_{\min}h_{2,\min}}\right\}.
\end{eqnarray}
\vskip2mm
\noindent
{\bf Theorem 3.1.} {\it Suppose that $G=(V,E)$ be a finite graph and the following conditions hold: \\
$(H)$\; $h_i(x)>0$ for all $x\in V$, $i=1,2$;\\
$(F_0)$\; $F(x,s,t)$ is continuously differentiable in $(s,t)\in \R^2$ for all $x\in V$;\\
$(F_1)$ \;   $\int_VF(x,0,0)d\mu=0$;\\
$(F_2)$ \;
$$0<A_V:=\liminf_{y\rightarrow +\infty}\frac{\int_V \max_{|s|+|t|\leq y}F(x,s,t)d\mu}{y^\delta}<
\limsup_{|s|+|t|\rightarrow \infty}\frac{\int_V F(x,s,t)d\mu}{|s|^p+|t|^q}:=B_V,$$
where $\delta=\min\{p,q\}$ .\\
Then for each $\lambda\in(\lambda_{1,V},\lambda_{2,V})$ with $\lambda_{1,V}=\frac{\varrho_V}{B_V}$ and $\lambda_{2,V}=\frac{1}{p2^{p-1}K_VA_V}$, system (\ref{eq2}) possesses an unbounded sequence of solutions.}

\vskip2mm
  \par
 In order to prove Theorem 3.1, we work in  the space $W_V:=W^{m_1,p}(V)\times W^{m_2,q}(V)$ equipped with the norm  $\|(u,v)\|_V=\|u\|_{W^{m_1,p}(V)}+\|v\|_{W^{m_2,q}(V)}$. Then $(W_V,\|\cdot\|_V)$ is a finite dimensional  Banach space.
  \par
  Consider the functional $I_{\lambda,V}:W_V\to\R$ as
\begin{eqnarray*}
\label{EQ1} I_{\lambda,V}(u,v)=\frac{1}{p}\int_V(|\nabla^{m_1}u|^p+h_1(x)|u|^p)d\mu+\frac{1}{q}\int_V(|\nabla^{m_2}v|^q+h_2(x)|v|^q)d\mu-\lambda\int_V F(x,u,v)d\mu.
\end{eqnarray*}
Then under the assumptions of Theorem 3.1, $I_{\lambda,V}\in C^1(W_V,\R)$ and
\begin{eqnarray*}\label{EQ2}
\langle I_{\lambda,V}'(u,v),(\phi_1,\phi_2)\rangle&=&\int_V\left[(\pounds_{m_1,p}u)\phi_1+h_1(x)|u|^{p-2}u\phi_1-\lambda  F_u(x,u,v)\phi_1\right]d\mu\nonumber\\
&&+\int_V\left[(\pounds_{m_2,q}v)\phi_2+h_2(x)|v|^{q-2}v\phi_2-\lambda F_v(x,u,v)\phi_2\right]d\mu
\end{eqnarray*}
for any $(u,v),(\phi_1,\phi_2)\in W_V$.  \par
A standard argument implies that $(u,v)\in W_V$ is a critical point of $I_{\lambda, V}$ iff
\begin{eqnarray*}
\int_V\left(\pounds_{m_1,p}u+h_1(x)|u|^{p-2}u-\lambda F_u(x,u,v)\right)\phi_1d\mu=0
\end{eqnarray*}
 and
\begin{eqnarray*}
\int_V\left(\pounds_{m_2,q}v+h_2(x)|v|^{q-2}v-\lambda F_v(x,u,v)\right)\phi_2d\mu=0
\end{eqnarray*}
for all $(\phi_1,\phi_2)\in W_V$. Furthermore,  by the arbitrary of $\phi_1$ and $\phi_2$, it can be  achieved that
\begin{eqnarray*}
\pounds_{m_1,p}u+h_1(x)|u|^{p-2}u=\lambda F_u(x,u,v),\\
\pounds_{m_2,q}v+h_2(x)|v|^{q-2}v=\lambda F_v(x,u,v).
\end{eqnarray*}
Therefore, seeking the solutions for system (\ref{eq2}) is equivalent to seeking the critical points of $I_{\lambda,V}$ on $W_V$ (see \cite{Zhang 2022} for example).
  \par
In order to apply Lemma 2.4, we shall exploit the functionals $\Phi_V:W_V\rightarrow\mathbb{R}$ and $\Psi_V:W_V\rightarrow\mathbb{R}$ which are set by
  \begin{eqnarray*}
\Phi_V(u,v)&=&\frac{1}{p}\int_V(|\nabla^{m_1}u|^p+h_1(x)|u|^p)d\mu+\frac{1}{q}\int_V(|\nabla^{m_2}v|^q+h_2(x)|v|^q)d\mu\\
&=&\frac{1}{p}\|u\|^p_{W^{m_1,p}(V)}+\frac{1}{q}\|v\|^q_{W^{m_2,q}(V)}
  \end{eqnarray*}
  and
\begin{eqnarray*}
\Psi_V(u,v)=\int_V F(x,u,v)d\mu.
\end{eqnarray*}
Then $I_{\lambda,V}(u,v)=\Phi_V-\lambda\Psi_V$. For each $r>\inf\limits_{W_V}\Phi_V$, define
\begin{eqnarray*}
\varphi_V(r)=\inf_{(u,v)\in\Phi_V^{-1}([-\infty,r])}\frac{\left(\sup_{(u,v)\in\Phi_V^{-1}([-\infty,r])}\Psi_V(u,v)\right)-\Psi_V(u,v)}{r-\Phi_V(u,v)}.
\end{eqnarray*}
  \par

 \vskip2mm
 \noindent
{\bf Lemma 3.1.}  {\it Assume that $(F_2)$ holds. Then $\gamma_V:=\liminf_{r\rightarrow+\infty}\varphi_V(r)<+\infty$.
}

 \vskip0mm
 \noindent
{\bf Proof.} Let $\{c_n\}$ be a real sequence satisfying $\lim_{n\rightarrow\infty}c_n=+\infty$ and
$$\lim_{n\rightarrow\infty}\frac{\int_V\max_{|s|+|t|\leq c_n}F(x,s,t)d\mu}{c_n^\delta}=A_V.$$
Write
$$r_n=\frac{2^{1-p}c_n^\delta}{pK_V},\;\;\mbox{for every}\;\; n\in \mathbb{N}.$$
By Lemma 2.1, for all $(u,v)\in W$ with $\Phi_V(u,v)\leq r_n$, we get
\begin{eqnarray}
\label{oo3.3}
     \frac{\|u\|_{\infty,V}^p}{p}+\frac{\|v\|_{\infty,V}^q}{q}
\leq K_V\left( \frac{\|u\|^p_{W^{m_1,p}(V)}}{p}+\frac{\|v\|^q_{W^{m_2,q}(V)}}{q}\right)
\leq K_V r_n,
\end{eqnarray}
\par
Next,  we claim that there exists $n_0\in \mathbb N$ such that $|u(x)|+|v(x)|\leq c_n$ for all $n\ge n_0$, all $x\in V$ and all $(u,v)\in W_V$ with $\Phi_V(u,v)\leq r_n$.  We prove the claim through the following three cases. Without loss of generality, we let $\delta=q$.
\par
(1) Assume that $\|u\|_{\infty,V}< 1$ and $\|v\|_{\infty,V}< 1$. It is obvious that there exists a $n_1\in \mathbb{N}$ such that $\|u\|_{\infty,V}+\|v\|_{\infty,V}\leq c_n$ for all $n>n_1$ by the fact $\lim_{n\rightarrow\infty}c_n=+\infty$.
\par
(2) Assume that $\|u\|_{\infty,V}\geq 1,\|v\|_{\infty,V}\geq 1$ or $\|u\|_{\infty,V}\geq 1,\|u\|_{\infty,V}< 1$. Then
\begin{eqnarray*}
     \frac{\|u\|_{\infty,V}^p}{p}+\frac{\|v\|_{\infty,V}^q}{q}
&  \geq  &  \frac{\|u\|_{\infty,V}^q+\|v\|_{\infty,V}^q}{p}\nonumber\\
&  \geq  &  \frac{2^{1-q}(\|u\|_{\infty,V}+\|v\|_{\infty,V})^q}{p},
\end{eqnarray*}
which together with (\ref{oo3.3}), implies that
\begin{eqnarray}
c_n^q\geq 2^{p-q}(\|u\|_{\infty,V}+\|v\|_{\infty,V})^q\geq (\|u\|_{\infty,V}+\|v\|_{\infty,V})^q.
\end{eqnarray}
Thus $\|u\|_{\infty,V}+\|v\|_{\infty,V}\leq c_n.$\\
\par
(3) Assume that $\|u\|_{\infty,V}< 1$ and $\|v\|_{\infty,V}\geq 1$. Then
\begin{eqnarray}
\label{3.5}
     \frac{\|u\|_{\infty,V}^p}{p}+\frac{\|v\|_{\infty,V}^q}{q}
&  \geq  &  \min\left\{\frac{1}{p},\frac{\|v\|_{\infty,V}^{q-p}}{q}\right\}(\|u\|_{\infty,V}^p+\|v\|_{\infty,V}^p)\nonumber\\
&  \geq  &  \min\left\{\frac{1}{p},\frac{\|v\|_{\infty,V}^{q-p}}{q}\right\}2^{1-p}(\|u\|_{\infty,V}+\|v\|_{\infty,V})^p.
\end{eqnarray}
\par
If $\min\{\frac{1}{p},\frac{\|v\|_{\infty,V}^{q-p}}{q}\}=\frac{1}{p}$, by (\ref{oo3.3}) and (\ref{3.5}), we have
\begin{eqnarray*}
    &&K_V r_n=\frac{2^{1-p}c_n^q}{p}\\
&  \geq  &  \frac{2^{1-p}(\|u\|_{\infty,V}+\|v\|_{\infty,V})^p}{p}\\
&  \geq  & \frac{2^{1-p}(\|u\|_{\infty,V}+\|v\|_{\infty,V})^q}{p}.
\end{eqnarray*}
Thus $\|u\|_{\infty,V}+\|v\|_{\infty,V}\leq c_n.$
\par
 By (\ref{oo3.3}), we have
\begin{eqnarray*}
\frac{\|v\|_{\infty,V}^q}{q}\leq \frac{2^{1-p}c_n^q}{p}.
\end{eqnarray*}
Note that $q-p\geq 0$. Then the above inequality implies that
\begin{eqnarray*}
\|v\|_{\infty,V}^{q-p} \geq \left(\frac{q}{p}\right)^\frac{q-p}{q}2^\frac{(1-p)(q-p)}{q}c_n^{q-p}.
\end{eqnarray*}
Thus, if $\min\{\frac{1}{p},\frac{\|v\|_{\infty,V}^{q-p}}{q}\}=\frac{\|v\|_{\infty,V}^{q-p}}{q}$, by (\ref{oo3.3}) and (\ref{3.5}), we have
\begin{eqnarray*}
K_V r_n& = &\frac{2^{1-p}c_n^q}{p}\nonumber\\
& \geq & \frac{\|v\|_{\infty,V}^{q-p}}{q}2^{1-p}(\|u\|_{\infty,V}+\|v\|_{\infty,V})^p\nonumber\\
& \geq & \left(\frac{q}{p}\right)^\frac{q-p}{q}2^\frac{(1-p)(q-p)}{q}c_n^{q-p}\frac{2^{1-p}}{q}(\|u\|_{\infty,V}+\|v\|_{\infty,V})^p.
\end{eqnarray*}
Hence, an easy computation implies that
\begin{eqnarray*}
c_n^p\geq \left(\frac{p}{q}\right)^\frac{p}{q}2^\frac{(1-p)(q-p)}{q}(\|u\|_{\infty,V}+\|v\|_{\infty,V})^p\geq (\|u\|_{\infty,V}+\|v\|_{\infty,V})^p.
\end{eqnarray*}
Thus based on the three cases, we conclude that for all $(u,v)\in W_V$ with $\Phi_V(u,v)\leq r_n$, we have $|u(x)|+|v(x)|\leq c_n$ for all $x\in V$.
Therefore, it follows from $(F_1)$ that
\begin{eqnarray*}
&   &   \varphi_V(r_n)\\
& = & \inf_{\frac{1}{p}\|u\|^p_{W^{m_1,p}(V)}+\frac{1}{q}\|v\|^q_{W^{m_2,q}(V)}\leq r_n}
\frac{\sup_{\frac{1}{p}\|u\|^p_{W^{m_1,p}(V)}+\frac{1}{q}\|v\|^q_{W^{m_2,q}(V)}\leq r_n}\int_VF(x,u,v)d\mu-\int_VF(x,u,v)d\mu}
{r_n-(\frac{1}{p}\|u\|^p_{W^{m_1,p}(V)}+\frac{1}{q}\|v\|^q_{W^{m_2,q}(V)})}\\
& \leq & \frac{\sup_{\frac{1}{p}\|u\|^p_{W^{m_1,p}(V)}+\frac{1}{q}\|v\|^q_{W^{m_2,q}(V)}\leq r_n}\int_VF(x,u,v)d\mu}{r_n}\\
&   =  & pK_V2^{p-1}\frac{\sup_{\frac{1}{p}\|u\|^p_{W^{m_1,p}(V)}+\frac{1}{q}\|v\|^q_{W^{m_2,q}(V)}\leq r_n}\int_VF(x,u,v)d\mu}{c_n^\delta}\\
& \leq & pK_V2^{p-1}\frac{\int_V\max_{|s|+|t|\leq c_n} F(x,s,t)d\mu}{c_n^\delta}.
\end{eqnarray*}
Hence,  $(F_2)$ implies that
$$\gamma_V \leq \liminf_{n\rightarrow\infty}\varphi_V(r_n) \leq  pK_V2^{p-1}A_V < pK_V2^{p-1}B_V \leq +\infty.$$
This finish the proof of the lemma. \qed

\vskip2mm
 \noindent
{\bf Lemma 3.2.} {\it For any fixed $\lambda\in(\lambda_{1,V},\lambda_{2,V})$,   $I_{\lambda,V}(u,v)=\Phi_V(u,v)-\lambda\Psi_V(u,v)$ is unbounded from below.}

\vskip0mm
 \noindent
{\bf Proof.} Assume that $\{\xi_n\}$ and $\{\eta_n\}$ are two positive real sequence satisfying $\lim_{n\rightarrow\infty}|\xi_n|+|\eta_n|=+\infty$ and
\begin{eqnarray}
\label{B_V}
\limsup_{n\rightarrow\infty}\frac{\int_V F(x,\xi_n,\eta_n)d\mu}{\xi_n^p+\eta_n^q}=B_V.
\end{eqnarray}
For each $n\in \mathbb{N}$, we  define
\begin{eqnarray*}
u_n(x)\equiv\xi_n,
v_n(x)\equiv\eta_n,\ \ \forall x\in V.
\end{eqnarray*}
It is obvious that  $(u_n,v_n)\in W_V$, $|\nabla^{m_i}u_n|=0$ and $|\nabla^{m_i}v_n|=0$ for all $m_i\ge 1,i=1,2.$
Then for every $n\in\mathbb{N}$, we have
\begin{eqnarray}\label{o1}
     I_{\lambda,V}(u_n,v_n)
& = &   \Phi_V(u_n,v_n)-\lambda\Psi_V(u_n,v_n)\nonumber\\
& = &   \frac{\xi_n^p}{p}\int_Vh_1(x)d\mu+\frac{\eta_n^q}{q}\int_Vh_2(x)d\mu-\lambda\int_V F(x,\xi_n,\eta_n)d\mu\nonumber\\
& \leq &  \varrho_V(\xi_n^p+\eta_n^q)-\lambda\int_V F(x,\xi_n,\eta_n)d\mu,
\end{eqnarray}
where $\varrho_V$ is defined by (\ref{3.1}).
\par
If $B_V<+\infty$, choosing $\epsilon_\lambda\in (\frac{\varrho_V}{\lambda B},1)$, by (\ref{B_V}), there exists $n_{\epsilon_\lambda}>0$ such that
$$\int_V F(x,\xi_n,\eta_n)d\mu > \epsilon_\lambda B_V(\xi_n^p+\eta_n^q), \;\;\;\forall n>n_{\epsilon_\lambda}.$$
Then combining with (\ref{o1}), we get
\begin{eqnarray*}
\Phi_V(u_n,v_n)-\lambda\Psi_V(u_n,v_n)
& \leq &  \varrho_V(\xi_n^p+\eta_n^q)-\lambda \epsilon_\lambda B_V(\xi_n^p+\eta_n^q)\\
&  =   &  (\varrho_V-\lambda \epsilon_\lambda B_V)(\xi_n^p+\eta_n^q), \ \ \ \forall n>n_{\epsilon_\lambda}.
\end{eqnarray*}
Then
$$\lim_{n\rightarrow\infty}[\Phi_V(u_n)-\lambda\Psi_V(v_n)]=-\infty.$$
\par
If $B_V=+\infty$, we consider $M_\lambda>\frac{\varrho_V}{\lambda}$. By (\ref{B_V}), there exists $n_{M_\lambda}$ such that
$$\int_V F(x,\xi_n,\eta_n)d\mu > {M_\lambda} (\xi_n^p+\eta_n^q),\;\;\;\forall n>n_{M_\lambda}.$$
Then combining with (\ref{o1}), we get
\begin{eqnarray*}
\Phi_V(u_n,v_n)-\lambda\Psi_V(u_n,v_n)
& \leq &  \varrho_V(\xi_n^p+\eta_n^q)-\lambda M_\lambda (\xi_n^p+\eta_n^q)\\
&  =   &  (\varrho_V-\lambda M_\lambda)(\xi_n^p+\eta_n^q), \ \ \ \forall n>n_{M_\lambda}.
\end{eqnarray*}
Noticing the choice of $M_\lambda$, we also have
$$\lim_{n\rightarrow\infty}[\Phi_V(u_n)-\lambda\Psi_V(v_n)]=-\infty.$$
Thus, we finish the proof of the lemma.\qed

\vskip2mm
\noindent
{\bf Lemma 3.3.} {\it $\Phi_V$ is sequentially weakly lower semi-continuous.}
\vskip0mm
\noindent
{\bf Proof.}  The proof is easily finished by exploiting the weak lower semi-continuity of the norm. \qed

\vskip2mm
\noindent
{\bf Lemma 3.4.} {\it $\Psi_V$ is sequentially weakly upper semi-continuous.}
\vskip0mm
\noindent
{\bf Proof.}  Assume that $(u_n,v_n)\rightharpoonup (u_0,v_0)$ in $W_V$. Note that $W_V$ is of finite dimension. Then $(u_n,v_n)\rightarrow (u_0,v_0)$ in $W_V$.
By $(F_0)$ and  the fact that $V$ is a finite set, it is easy to obtain that
\begin{eqnarray*}
      \lim_{n\rightarrow\infty}\int_VF(x,u_k,v_k)d\mu
&  =  & \lim_{n\rightarrow\infty}\sum_{x\in V}F(x,u_k,v_k)\mu(x)\\
&  =  &  \sum_{x\in V}F(x,u_0,v_0)\mu(x)\\
&  =  &  \int_VF(x,u_0,v_0)d\mu.
\end{eqnarray*}
Hence, $\Psi_V$ is  sequentially weakly upper semi-continuous in $W_V$.
\qed
\vskip2mm
\noindent
{\bf Proof of Theorem 3.1.} It is easy to see that $\Phi_V:W_V\to \R$ is coercive.  Lemma 3.1--Lemma 3.4 imply that all of conditions in Lemma 2.4 are satisfied. Hence, Lemma 2.4 (a) implies that for each $\left( \lambda_{1,V}  \lambda_{2,V} \right)$,  the functional $I_{\lambda,V}$ has a sequence $\{(u_n^*,v_n^*)\}$ of critical points that are solutions of system (\ref{eq2}) such that $\lim_{n\to\infty}\Phi_V(u_n^*,v_n^*)=+\infty$.
\qed

\vskip2mm
{\section{Result and proofs for  system (\ref{eq3})}}
  \setcounter{equation}{0}
In this section, we investigate the generalized poly-Laplacian system (\ref{eq3}) and obtain the following result.
\par
Let
$$K_\Omega=\max\left\{\frac{C^p(m_1,p,\Omega)}{\mu_{\min,\Omega}},\frac{C^q(m_2,q,\Omega)}{\mu_{\min,\Omega}}\right\},$$
where $C(m_1,p,\Omega$ and $C(m_2,q,\Omega)$ are defined in Lemma 2.2.
\vskip2mm
  \noindent
{\bf Theorem 4.1.} {\it Assume that $G=(V,E)$ is a locally finite graph and the following conditions hold: \\
 $(H)'$\: $h_i(x)>0$ for all $x\in \Omega$, $i=1,2$;\\
 $(F_0)'$\; $F(x,s,t)$ is continuously differentiable in $(s,t)\in \R^2$ for all $x\in \Omega$;\\
$(F_1)'$ \;   $\int_{\Omega}F(x,0,0)d\mu=0$;\\
$(F_2)'$ \;
$$ 0< A_\Omega:= \liminf_{y\rightarrow +\infty}\frac{\int_{\Omega}\max_{|s|+|t|\leq y}F(x,s,t)d\mu}{y^\delta}<
\limsup_{|s|+|t|\rightarrow\infty}\frac{\int_{\Omega} F(x,s,t)d\mu}{|s|^p+|t|^q} :=B_\Omega.$$
where $\delta=\min\{p,q\}$ .\\
Then for each $\lambda\in (\lambda_{1,\Omega},\lambda_{2,\Omega})$ with $\lambda_{1,\Omega}=\frac{1}{B_\Omega}$ and $\lambda_{2,\Omega}=\frac{1}{pK_\Omega2^{p-1}A_\Omega}$, system (\ref{eq3}) possesses an unbounded sequence of solutions.}

 \vskip2mm
 \par
 The proofs of Theorem 4.1 is the essentially same as Theorem 3.1 with some slight modifications. In order to prove Theorem 4.1, we work in  the space $W_0:=W_0^{m_1,p}(\Omega)\times W_0^{m_2,q}(\Omega)$ equipped with the norm   $\|(u,v)\|_0=\|u\|_{W_0^{m_1,p}(\Omega)}+\|v\|_{W_0^{m_2,q}(\Omega)}$. Then $(W_0,\|\cdot\|_0)$ is a  finite dimensional  Banach space. Consider the functional $I_{\lambda,\Omega}:W_0\to\R$ as
\begin{eqnarray}
\label{4.1} I_{\lambda,\Omega}(u,v)=\frac{1}{p}\int_{\Omega\cup\partial\Omega}|\nabla^{m_1}u|^pd\mu+\frac{1}{q}\int_{\Omega\cup\partial\Omega}|\nabla^{m_2}v|^qd\mu
-\lambda\int_\Omega F(x,u,v)d\mu.
\end{eqnarray}
Then under the assumptions of Theorem 4.1, $I_{\lambda,\Omega}\in C^1(W_0,\R)$ and
\begin{eqnarray}\label{4.2}
          \langle I_{\lambda,\Omega}'(u,v),(\phi_1,\phi_2)\rangle
&  =  &  \int_{\Omega\cup\partial\Omega}\left[(\pounds_{m_1,p}u)\phi_1-\lambda F_u(x,u,v)\phi_1\right]d\mu\nonumber\\
&     &  +\int_{\Omega\cup\partial\Omega}\left[(\pounds_{m_2,q}v)\phi_2-\lambda F_v(x,u,v)\phi_2\right]d\mu
\end{eqnarray}
for any $(u,v),(\phi_1,\phi_2)\in W_0$.
\par
Obviously, $(u,v)\in W_0$ is a critical point of $I_{\lambda,\Omega}$ iff
\begin{eqnarray*}
\int_{\Omega\cup\partial\Omega}\left(\pounds_{m_1,p}u-\lambda F_u(x,u,v)\right)\phi_1d\mu=0
\end{eqnarray*}
 and
\begin{eqnarray*}
\int_{\Omega\cup\partial\Omega}\left(\pounds_{m_2,q}v-\lambda F_v(x,u,v)\right)\phi_2d\mu=0
\end{eqnarray*}
for all $(\phi_1,\phi_2)\in W_0$. Furthermore,  by the arbitrary of $\phi_1$ and $\phi_2$, it can be  achieved that
system (\ref{eq3}) holds.
Therefore, seeking the solutions for system (\ref{eq3}) is equivalent to seeking the critical points of $I_{\lambda,V}$ on $W_0$.
\par
In order to apply Lemma 2.4, we will use the functionals $\Phi_\Omega:W_0\rightarrow\mathbb{R}$ and $\Psi_\Omega:W_0\rightarrow\mathbb{R}$ defined by
  \begin{eqnarray*}
\Phi_\Omega(u,v)&=&\frac{1}{p}\int_{\Omega\cup\partial\Omega}|\nabla^{m_1}u|^pd\mu+\frac{1}{q}\int_{\Omega\cup\partial\Omega}|\nabla^{m_2}v|^qd\mu\\
&=&\frac{1}{p}\|u\|^p_{W_0^{m_1,p}(\Omega)}+\frac{1}{q}\|v\|^q_{W_0^{m_2,q}(\Omega)}
  \end{eqnarray*}
  and
\begin{eqnarray*}
\Psi_\Omega(u,v)=\int_\Omega F(x,u,v)d\mu.
\end{eqnarray*}
Then $I_{\lambda,\Omega}(u,v)=\Phi_\Omega-\lambda\Psi_\Omega$ and for every $r>\inf\limits_{W_0}\Phi_\Omega$, define
\begin{eqnarray*}
\varphi_\Omega(r)=\inf_{(u,v)\in\Phi_\Omega^{-1}([-\infty,r])}\frac{\left(\sup_{(u,v)\in\Phi_\Omega^{-1}([-\infty,r])}
\Psi_\Omega(u,v)\right)-\Psi_\Omega(u,v)}{r-\Phi_\Omega(u,v)}.
\end{eqnarray*}

\vskip2mm
 \noindent
{\bf Lemma 4.1.}  {\it Assume that $(F_2)'$ holds. Then $\gamma_\Omega:=\liminf_{r\rightarrow+\infty}\varphi_\Omega(r)<+\infty$.
}

 \vskip0mm
 \noindent
{\bf Proof.} The proof is the same as Theorem 3.1 with substituting $\Omega$, $K_\Omega$, $A_\Omega$, $B_\Omega$, $\|u\|_{\infty, \Omega}$ and $\|v\|_{\infty, \Omega}$  for $V$, $K_V$, $A_V$, $B_V$, $\|u\|_{\infty, V}$ and $\|v\|_{\infty, V}$, respectively. We omit the details. \qed

\vskip2mm
 \noindent
{\bf Lemma 4.2.} {\it For any fixed $\lambda\in(\lambda_{1,\Omega},\lambda_{2,\Omega})$, $I_{\lambda,\Omega}(u,v)=\Phi_\Omega(u,v)-\lambda\Psi_\Omega(u,v)$ is unbounded from below.}

\vskip0mm
 \noindent
{\bf Proof.} Suppose that $\{\xi_n\}$ and $\{\eta_n\}$ are two positive real sequence such that $\lim_{n\rightarrow\infty}|\xi_n|+|\eta_n|=+\infty$ and
\begin{eqnarray}
\label{B_o}
\lim_{n\rightarrow\infty}\frac{\int_\Omega F(x,\xi_n,\eta_n)d\mu}{\xi_n^p+\eta_n^q}=B_\Omega.
\end{eqnarray}
For each $n\in \mathbb{N}$, we define
\begin{eqnarray*}
u_n(x)\equiv\xi_n,
v_n(x)\equiv\eta_n,\ \ \forall x\in \Omega.
\end{eqnarray*}
It is easy to check that  $(u_n,v_n)\in W_0$, $|\nabla^{m_i}u_n|=0$ and $|\nabla^{m_i}v_n|=0$ for all $m_i\ge 1,i=1,2.$
Then
\begin{eqnarray*}
     I_{\lambda,\Omega}(u_n,v_n)
& = &   \Phi_\Omega(u_n,v_n)-\lambda \Psi_\Omega(u_n,v_n)\\
& = &   -\int_\Omega F(x,\xi_n,\eta_n)d\mu.
\end{eqnarray*}
\par
If $B_\Omega<+\infty$, choosing  $\epsilon_\lambda\in (0,\frac{1}{\lambda B_\Omega})$, by (\ref{B_o}), there exists $n_{\epsilon_\lambda}$ such that
$$\int_\Omega F(x,\xi_n,\eta_n)d\mu > \epsilon_\lambda B_\Omega(\xi_n^p+\eta_n^q), \;\;\;\forall n>n_{\epsilon_\lambda}.$$
Hence
\begin{eqnarray*}
\Phi_\Omega(u_n,v_n)-\lambda\Psi_\Omega(u_n,v_n)
& \leq &  -\lambda\epsilon_\lambda B_\Omega(\xi_n^p+\eta_n^q), \ \ \ \forall n>n_{\epsilon_\lambda}.
\end{eqnarray*}
Thus,
$$\lim_{n\rightarrow\infty}[ \Phi_\Omega(u_n)-\lambda \Psi_\Omega(v_n)]=-\infty.$$
\par
If $B_\Omega=+\infty$,  consider $M_\lambda>\frac{1}{\lambda}$. By (\ref{B_o}), there exists $n_{M_\lambda}$ such that
$$\int_\Omega F(x,\xi_n,\eta_n)d\mu > M_\lambda (\xi_n^p+\eta_n^q),\;\;\;\forall n>n_{M_\lambda}.$$
Hence
\begin{eqnarray*}
 \Phi_\Omega(u_n,v_n)-\lambda \Psi_\Omega(u_n,v_n)
& \leq &  -\lambda M_\lambda (\xi_n^p+\eta_n^q), \ \ \ \forall n>n_{M_\lambda}.
\end{eqnarray*}
By the choice of $M_\lambda$, we also have
$$\lim_{n\rightarrow\infty}[ \Phi_\Omega(u_n)-\lambda\Psi_\Omega(v_n)]=-\infty.$$
Thus, we finish the proof of this lemma. \qed

\vskip2mm
\noindent
{\bf Lemma 4.3.} {\it $ \Phi_\Omega$ is sequentially weakly lower semi-continuous.}
\vskip0mm
\noindent
{\bf Proof.}  The proof is easily completed by using the weak lower semi-continuity of the norm. \qed

\vskip2mm
\noindent
{\bf Lemma 4.4.} {\it $ \Psi_\Omega$ is sequentially weakly upper semi-continuous.}
\vskip0mm
\noindent
{\bf Proof.}  The proof is the same as Lemma 3.4 with replacing $W$ with $W_0$ and $V$ with $\Omega$.\qed

\vskip2mm
\noindent
{\bf Proof of Theorem 4.1.} It is obvious that $ \Phi_\Omega:W_0\to \R$ is coercive.
 Lemma 4.1--Lemma 4.4 imply that all of  conditions in Lemma 2.4 are satisfied. Hence, Lemma 2.4(a) implies that for each $\left( \lambda_{1,\Omega},\lambda_{2,\Omega}\right)$,  $I_{\lambda,\Omega}$ has a sequence $\{(u_n^{\star},v_n^{\star})\}$ of critical points that are solutions of system (\ref{eq3}) such that $\lim_{n\to\infty}\Phi_\Omega(u_n^{\star},v_n^{\star})=+\infty$.
\qed

\vskip2mm
{\section{Result and proofs for system (\ref{eq4})}}
  \setcounter{equation}{0}
  In this section, we investigate the $(p,q)$-Laplacian system (\ref{eq4}). We first make the following assumptions:\\
  $(M_1)$\: there exists a $\mu_0>0$ such that $\mu(x)\ge \mu_0$ for all $x\in V$;\\
  $(M_2)$ there exists a $x_0\in V$ such that $M_1(x_0)\leq M_1(x)$ and $M_2(x_0)\leq M_2(x)$ for all $x\in V$, where
  \begin{eqnarray*}
M_1(x)& = & \left(\frac{deg(x)}{2\mu(x)}\right)^{\frac{p}{2}}\mu(x)+h_1(x)\mu(x)+\sum_{y\thicksim x}\left(\frac{w_{xy}}{2\mu(y)}\right)^{\frac{p}{2}}\mu(y),\ \ x\in V,\nonumber\\
M_2(x)& = & \left(\frac{deg(x)}{2\mu(x)}\right)^{\frac{q}{2}}\mu(x)+h_2(x)\mu(x)+\sum_{y\thicksim x}\left(\frac{w_{xy}}{2\mu(y)}\right)^{\frac{q}{2}}\mu(y),\ \ x\in V.
\end{eqnarray*}
  $(H_1)$\; there exists a constant $h_0>0$ such that $h_i(x)\geq h_0>0$ for all $x\in V$, $i=1,2$;\\
  \par
  Let
  \begin{eqnarray} \label{oo1}
   \varrho=\max\left\{\frac{M_1(x_0)}{p},\frac{M_2(x_0)}{q}\right\}\;\;\mbox{and}\ \ K=\max\left\{\frac{1}{h_0^{1/p}\mu_0^{1/p}},\frac{1}{h_0^{1/q}\mu_0^{1/q}}\right\}.
  \end{eqnarray}
  \vskip2mm
  \noindent
{\bf Theorem 5.1.} {\it Suppose that $G=(V,E)$ is a locally finite graph, and  $(M_1)$,  $(M_2)$, $(H_1)$  and  the following conditions hold: \\
$(\tilde F_0)$\; $F(x,s,t)$ is continuously differentiable in $(s,t)\in \R^2$ for all $x\in V$,  and there exists a  function $a\in C(\R^+,\R^+)$ and a function $b:V\to \R^+$ with $b\in L^1(V)$ such that
$$
|F_s(x,s,t)|\le a(|(s,t)|) b(x), |F_t(x,s,t)|\le a(|(s,t)|) b(x), |F(x,s,t)|\le a(|(s,t)|) b(x),
$$
for all $x\in V$ and all $(s,t)\in \R^2$;\\
 $(\tilde F_1)$\; $\int_V F(x,0,0)d\mu=0$ ;\\
 $(\tilde F_2)$ \;
$$0<A:=\liminf_{y\rightarrow\infty}\frac{\int_V \max_{|s|+|t|\leq y}F(x,s,t)d\mu}{y^\delta}<
\limsup_{|s|+|t|\rightarrow +\infty}\frac{\int_V F(x,s,t)d\mu}{|s|^p+|t|^q}:=B.$$
where $\delta=\min\{p,q\}$.\\
Then for each  $\lambda\in(\Theta_1,\Theta_2)$ with  $\Theta_1=\frac{ \varrho}{B}$ and $\Theta_2=\frac{1}{pK2^{p-1}A}$, system (\ref{eq4})  possesses an unbounded sequence of solutions.}

\vskip2mm
\par
 We work in the space $W:=W_{h_1}^{1,p}(V)\times W_{h_2}^{1,q}(V)$ with the norm  equipped with $\|(u,v)\|=\|u\|_{W_{h_1}^{1,p}(V)}+\|v\|_{W_{h_2}^{1,q}(V)}$ and then $(W,\|\cdot\|)$ is a  Banach space which is infinite dimensional.
\par
We consider the functional $I_\lambda:W\to\R$ as
\begin{eqnarray}
\label{EQ1}
I_\lambda(u,v)=\frac{1}{p}\int_V(|\nabla u|^p+h_1(x)|u|^p)d\mu+\frac{1}{q}\int_V(|\nabla v|^q+h_2(x)|v|^q)d\mu-\lambda\int_V F(x,u,v)d\mu.
\end{eqnarray}
Then by Appendix A.2 in \cite{Yang 2023}, under the assumptions of Theorem 5.1,  $I_\lambda\in C^1(W,\mathbb{R})$, and
\begin{eqnarray}\label{EQ2}
\langle I_\lambda'(u,v),(\phi_1,\phi_2)\rangle&=&\int_V\left[|\nabla u|^{p-2}\Gamma(u,\phi_1)+h_1(x)|u|^{p-2}u\phi_1-\lambda  F_u(x,u,v)\phi_1\right]d\mu\nonumber\\
&&+\int_V\left[|\nabla v|^{q-2}\Gamma(v,\phi_2)+h_2(x)|v|^{q-2}v\phi_2-\lambda F_v(x,u,v)\phi_2\right]d\mu
\end{eqnarray}
for any $(u,v),(\phi_1,\phi_2)\in W$.
\par
 Obviously, $(u,v)\in W$ is a critical point of $I_\lambda$ iff
\begin{eqnarray*}
\int_V\left[|\nabla u|^{p-2}\Gamma(u,\phi_1)+h_1(x)|u|^{p-2}u\phi_1-\lambda F_u(x,u,v)\phi_1\right]d\mu=0
\end{eqnarray*}
 and
\begin{eqnarray*}
\int_V\left[|\nabla v|^{q-2}\Gamma(v,\phi_2)+h_2(x)|v|^{q-2}v\phi_2-\lambda F_v(x,u,v)\phi_2\right]d\mu=0
\end{eqnarray*}
for all $(\phi_1,\phi_2)\in W$. Furthermore,  by the arbitrary of $\phi_1$ and $\phi_2$, it can be  achieved that
system (\ref{eq4}) holds.
Therefore, seeking the solutions for system (\ref{eq4}) is equivalent to seeking the critical points of $I_{\lambda}$ on $W$.
\par
Define $ \Phi:W\rightarrow\mathbb{R}$ and $ \Psi:W\rightarrow\mathbb{R}$ by
  \begin{eqnarray*}
\Phi(u,v)&=&\frac{1}{p}\int_V(|\nabla u|^p+h_1(x)|u|^p)d\mu+\frac{1}{q}\int_V(|\nabla v|^q+h_2(x)|v|^q)d\mu\\
&=&\frac{1}{p}\|u\|^p_{W_h^{1,p}(V)}+\frac{1}{q}\|v\|^q_{W_h^{1,q}(V)}
  \end{eqnarray*}
  and
\begin{eqnarray*}
 \Psi(u,v)=\int_V F(x,u,v)d\mu.
\end{eqnarray*}
Then $I_\lambda(u,v)= \Phi-\lambda \Psi$. For every $r>\inf\Phi$, set
\begin{eqnarray*}
\varphi(r)=\inf_{(u,v)\in\Phi^{-1}([-\infty,r])}\frac{\left(\sup_{(u,v)\in\Phi^{-1}([-\infty,r])}\Psi(u,v)\right)-\Psi(u,v)}{r-\Phi(u,v)}.
\end{eqnarray*}

\vskip2mm
 \noindent
{\bf Lemma 5.1.}  {\it Assume that $(\tilde F_2)$ holds. Then $\gamma:=\liminf_{r\rightarrow+\infty}\varphi(r)<+\infty$.
}
\vskip0mm
 \noindent
 {\bf Proof.} The proof is the essentially same as Theorem 3.1 with substituting the locally finite graph $V$, $K$, $A$, $B$, $\|u\|_{\infty}$ and $\|v\|_{\infty}$  for finite graph $V$, $K_V$, $A_V$, $B_V$, $\|u\|_{\infty, V}$ and $\|v\|_{\infty, V}$, respectively. We omit the details. \qed

 \vskip2mm
 \noindent
{\bf Lemma 5.2.} {\it For any given $\lambda\in(\Theta_1,\Theta_2)$, the functional
 $I_\lambda(u ,v )=\Phi(u ,v )-\lambda \Psi(u ,v )$ is unbounded from below.}

\vskip0mm
 \noindent
{\bf Proof.} By $(\tilde F_2)$, we can assume that $\{\xi_n\}$ and $\{\eta_n\}$ are two positive real sequence satisfying $\lim_{n\rightarrow\infty}|\xi_n|+|\eta_n|=+\infty$ and
\begin{eqnarray}
\label{B_p}
\lim_{n\rightarrow\infty}\frac{\int_V F(x,\xi_n,\eta_n)d\mu}{\xi_n^p+\eta_n^q}=B.
\end{eqnarray}
For each $n\in \mathbb{N}$ define
\begin{eqnarray*}
u_n(x)=\begin{cases}
                  \xi_n,& x=x_0\\
                  0,&x\not= x_0
                  \end{cases},
\quad v_n(x)=\begin{cases}
                  \eta_n,& x=x_0\\
                  0,&x\not= x_0,
                  \end{cases}
\end{eqnarray*}
 where $x_0$ is given in the assumption $(M_2)$. Then a simple calculation implies that
$$
|\nabla u_n|(x)=\begin{cases}
           \sqrt{\frac{deg(x_0)}{2\mu(x_0)}}\xi_n, & x=x_0,\\
           \sqrt{\frac{w_{x_0y}}{2\mu(y)}}\xi_n, & x=y \mbox{ with } y\thicksim x_0,\\
           0, &\mbox{otherwise},
           \end{cases}
$$
and
$$
|\nabla v_n|(x)=\begin{cases}
           \sqrt{\frac{deg(x_0)}{2\mu(x_0)}}\eta_n, & x=x_0,\\
           \sqrt{\frac{w_{x_0y}}{2\mu(y)}}\eta_n, & x=y \mbox{ with } y\thicksim x_0,\\
           0, &\mbox{otherwise}.
           \end{cases}
$$
Then
\begin{eqnarray*}\label{pp1}
&     &     \int_V{(|\nabla u_n|^p+h_1(x)|u_n|^p)}d\mu\nonumber\\
&  =  &    \sum_{x\in V}(|\nabla u_n(x)|^p+h_1(x)|u_n(x)|^p)\mu(x)\nonumber\\
&  =  &    (|\nabla u_n(x_0|^p)+h_1(x_0)|u_n(x_0)|^p)\mu(x_0)+\sum_{y\thicksim x_0}(|\nabla u_n(y)|^p+ h_1(y)|u_n(y)|^p)\mu(y)\nonumber\\
&  =  & \left(\frac{deg(x_0)}{2\mu(x_0)}\right)^{\frac{p}{2}}\xi_n^p\mu(x_0)+h_1(x_0)\xi_n^p\mu(x_0)+\xi_n^p\sum_{y\thicksim x_0}\left(\frac{w_{x_0y}}{2\mu(y)}\right)^{\frac{p}{2}}\mu(y)\nonumber\\
&  =  & \xi_n^p M_1(x_0) ,
\end{eqnarray*}
and similarly,
\begin{eqnarray*}\label{pp2}
&     &     \int_V{(|\nabla v_n|^q+h_2(x)|v_n|^q)}d\mu\nonumber\\
&  =  & \left(\frac{deg(x_0)}{2\mu(x_0)}\right)^{\frac{q}{2}}\eta_n^q\mu(x_0)+h_2(x_0)\eta_n^q\mu(x_0)+\eta_n^q\sum_{y\thicksim x_0}\left(\frac{w_{x_0y}}{2\mu(y)}\right)^{\frac{q}{2}}\mu(y)\nonumber\\
&  =  & \eta_n^qM_2(x_0) ,
\end{eqnarray*}
where $M_1(x_0)$ and $M_2(x_0)$ are given in the assumption $(M_2)$.
Then $\{(u_n,v_n)\}\subset W$ and for every $n\in\mathbb{N}$, we have
\begin{eqnarray}
\label{5.11}
     I_\lambda(u_n,v_n)
& = &   \Phi(u_n,v_n)-\lambda \Psi(u_n,v_n)\nonumber\\
& = &   \frac{\xi_n^pM_1(x_0)}{p}+\frac{\eta_n^qM_2(x_0)}{q}-\int_V F(x,\xi_n,\eta_n)d\mu\nonumber\\
& \leq &   \varrho(\xi_n^p+\eta_n^q)-\int_V F(x,\xi_n,\eta_n)d\mu,
\end{eqnarray}
where $\varrho$ is given in (\ref{oo1}).
\par
If $B<+\infty$, choosing $\tilde\epsilon_\lambda \in (\frac{ \varrho}{\lambda B},1)$, by (\ref{B_p}), there exists $n_{\tilde\epsilon_\lambda}$ such that
$$\int_V F(x,\xi_n,\eta_n)d\mu > \tilde\epsilon_\lambda B(\xi_n^p+\eta_n^q), \;\;\;\forall n>n_{\tilde\epsilon_\lambda}.$$
Thus, combining with (\ref{5.11}), we have
\begin{eqnarray*}
 \Phi(u_n,v_n)-\lambda \Psi(u_n,v_n)
& \leq &   \varrho(\xi_n^p+\eta_n^q)-\lambda \tilde\epsilon_\lambda B(\xi_n^p+\eta_n^q)\\
&  =   &  ( \varrho-\lambda \tilde\epsilon_\lambda B)(\xi_n^p+\eta_n^q), \ \ \ \forall n>n_{\tilde\epsilon_\lambda}.
\end{eqnarray*}
 Hence
$$\lim_{n\rightarrow\infty}[\Phi(u_n)-\lambda\Psi(v_n)]=-\infty.$$
\par
If $B=+\infty$, let us consider $\tilde M_\lambda>\frac{ \varrho}{\lambda}$. By (\ref{B_p}), there exists $n_{\tilde M_\lambda}$ such that
$$\int_V F(x,\xi_n,\eta_n)d\mu > \tilde M_\lambda (\xi_n^p+\eta_n^q),\;\;\;\forall n>n_{\tilde M_\lambda}.$$
Thus
\begin{eqnarray*}
 \Phi(u_n,v_n)-\lambda \Psi(u_n,v_n)
& \leq &  \varrho(\xi_n^p+\eta_n^q)-\lambda \tilde M_\lambda (\xi_n^p+\eta_n^q)\\
&  =   &  ( \varrho-\lambda \tilde M_\lambda)(\xi_n^p+\eta_n^q), \ \ \ \forall n>n_{\tilde M_\lambda}.
\end{eqnarray*}
Combining the choice of $\tilde M_\lambda$, in this case, we also have
$$\lim_{n\rightarrow\infty}[ \Phi(u_n)-\lambda \Psi(v_n)]=-\infty.$$
Thus we complete the proof of this lemma.\qed

\vskip2mm
\noindent
{\bf Lemma 5.3.} {\it $ \Phi$ is sequentially weakly lower semi-continuous.}
\vskip0mm
\noindent
{\bf Proof.}  The proof is easily completed by using the weak lower semi-continuity of the norm. \qed

\vskip2mm
\noindent
{\bf Lemma 5.4.} {\it $ \Psi$ is sequentially weakly upper semi-continuous.}
\vskip0mm
\noindent
{\bf Proof.}  Assume that  $(u_k,v_k)\rightharpoonup (u_0,v_0)$ for some $(u_0,v_0)\in W$.  Then,
$$\lim_{k\rightarrow\infty}\int_Vu_k\varphi d\mu=\int_Vu_0\varphi d\mu,\forall\varphi\in C_c(V),$$
which implies that
\begin{eqnarray}
\lim_{k\rightarrow\infty}u_k(x)=u_0(x)\;\; \mbox{for any fixed}\;x\in V
\end{eqnarray}
by choosing
$$
\varphi(y)=
\begin{cases}
1,& y=x\\
0, & y\not=x.
\end{cases}
$$
Similarly, we have
\begin{eqnarray*}\label{eq14}
\lim_{k\rightarrow\infty}v_k(x)=v_0(x)\;\; \mbox{for any fixed}\;x\in V.
\end{eqnarray*}
By $(\tilde F_0)$ and Lebesgue dominated convergence theorem, it is easy to obtain that
\begin{eqnarray*}
      \lim_{n\rightarrow\infty}\int_VF(x,u_k,v_k)d\mu
  =    \int_VF(x,u_0,v_0)d\mu.
\end{eqnarray*}
Hence, $ \Psi$ is sequentially weakly upper semi-continuous in $W$.
\qed
\vskip2mm
\noindent
{\bf Proof of Theorem 5.1.} Obviously, $ \Phi:W \to \R$ is coercive.  Lemma 5.1--Lemma 5.4 imply that all of conditions in Lemma 2.4(a) hold for $I_\lambda$. Hence, Lemma 2.4(b) implies that for each  $\lambda\in\left( \Theta_1,  \Theta_2 \right)$, $I_\lambda$ has a sequence $\{(u_n,v_n)\}$ of critical points that are solutions of system (\ref{eq4}) such that $\lim_{n\to\infty}\Phi(u_n,v_n)=+\infty$.
\qed
\vskip2mm
{\section{The results for the scalar equations}}
 \setcounter{equation}{0}
\vskip2mm
\par
By using the similar arguments of Theorem 3.1, we can also obtain the similar results for the following scalar equation on finite graph $G=(V,E)$:
\begin{eqnarray}
\label{eqq21}
  \pounds_{m,p}u+h(x)|u|^{p-2}u=\lambda f(x,u),\;\;\;\;\hfill x\in V,
\end{eqnarray}
where $m\geq 1$ is an integer, $h:V\rightarrow\mathbb{R}$, $p>1$, $\lambda>0$ and  $f:V\times\mathbb{R}\rightarrow\mathbb{R}$.

\vskip2mm
\noindent
{\bf Theorem 6.1.} {\it Let $G=(V,E)$ be a finite graph and $F(x,s)=\int_0^s f(x,\tau)d\tau$ for all $x\in V$. Assume that the following conditions hold: \\
$(h)$\; $h(x)>0$ for all $x\in V$;\\
$(f_0)$\; $F(x,s)$ is continuously differentiable in $s\in \R$ for all $x\in V$;\\
$(f_1)$ \;   $\int_VF(x,0)d\mu=0$;\\
$(f_2)$ \;
$$ 0<\tilde {A}_V:=\liminf_{y\rightarrow\infty}\frac{\int_V \max_{|s|\leq y}F(x,s,t)d\mu}{y^p}<
 \limsup_{|s|\rightarrow\infty}\frac{\int_V F(x,s)d\mu}{|s|^p}:=\tilde {B}_V.$$
Then for each  $\lambda\in(\tilde\lambda_{1,V},\tilde\lambda_{2,V})$ with $\tilde\lambda_{1,V}=\frac{\tilde{\varrho}_V}{\tilde{B}_V}$ and $\tilde\lambda_{2,V}=\frac{1}{p\tilde {K}_V\tilde{A}_V}$, where
$\tilde {K}_V=\frac{1}{h_{\min}\mu_{\min}}$ and $\tilde{\varrho}_V=\frac{1}{p}\int_V h(x)d\mu$, equation (\ref{eqq21}) possesses an unbounded sequence of solutions.}

\vskip2mm
\par
The proofs of Theorem 6.1 are almost same as those of Theorem 3.1 and even more simple because there is  no the couple term. Here, we just present the proof that $\tilde{\gamma}_V:=\liminf_{r\rightarrow+\infty}\tilde{\varphi}_V(r)<+\infty$ which is related the range of the parameter of $\lambda$ and also show that the proof for single equation is indeed more simple, where
\begin{eqnarray*}
  \tilde{\varphi}_V(r)
  & = & \inf_{u\in\tilde{\Phi}_V^{-1}([-\infty,r])}\frac{\left(\sup_{(u,v)\in\tilde{\Phi}_V^{-1}([-\infty,r])}\tilde{\Psi}_V(u)\right)-\tilde{\Psi}_V(u)}{r-\tilde{\Phi}_V(u)},\\
\tilde{\Phi}_V(u) & = & \frac{1}{p}\int_V(|\nabla^{m}u|^p+h(x)|u|^p)d\mu=\frac{1}{p}\|u\|^p_{W^{m,p}(V)},\\
\tilde{\Psi}_V(u) & = &\int_VF(x,u)d\mu,
\end{eqnarray*}
and $\tilde{I}_{\lambda,V}=\tilde{\Phi}_V-\lambda\tilde{\Psi}_V$ is the corresponding variational functional of (\ref{eqq21}).
\par
In fact, let $\{c_n\}$ be a real sequence satisfying $\lim_{n\rightarrow\infty}c_n=+\infty$ and
$$\lim_{n\rightarrow\infty}\frac{\int_V\max_{|s|\leq c_n}F(x,s)d\mu}{c_n^\delta}=\tilde{A}_V.$$
Write
$$r_n=\frac{c_n^p}{p\tilde{K}_V},\;\;\mbox{for every}\;\; n\in \mathbb{N}.$$
By Lemma 2.1, for all $(u,v)\in W$ with $\tilde{\Phi}_V(u)\leq r_n$, we get
\begin{eqnarray*}
     \frac{\|u\|_{\infty,V}^p}{p}
\leq \tilde{K}_V\frac{\|u\|^p_{W^{m,p}(V)}}{p}
\leq \tilde{K}_V r_n,
\end{eqnarray*}
Hence, $|u(x)|\le c_n$ for all $x\in V$. Therefore, it follows from $(f_1)$ that
\begin{eqnarray*}
&   &   \tilde{\varphi}_V(r_n)\\
& = & \inf_{\frac{1}{p}\|u\|^p_{W^{m,p}(V)}\leq r_n}
\frac{\sup_{\frac{1}{p}\|u\|^p_{W^{m_1,p}(V)}\leq r_n}\int_VF(x,u)d\mu-\int_VF(x,u)d\mu}
{r_n-\frac{1}{p}\|u\|^p_{W^{m_1,p}(V)}}\\
& \leq & \frac{\sup_{\frac{1}{p}\|u\|^p_{W^{m_1,p}(V)}\leq r_n}\int_VF(x,u)d\mu}{r_n}\\
&   =  & p\tilde{K}_V\frac{\sup_{\frac{1}{p}\|u\|^p_{W^{m_1,p}(V)}\leq r_n}\int_VF(x,u)d\mu}{c_n^\delta}\\
& \leq & p\tilde{K}_V\frac{\int_V\max_{|s|\leq c_n} F(x,s)d\mu}{c_n^\delta}.
\end{eqnarray*}
Hence,  $(f_2)$ implies that
$$\tilde{\gamma}_V \leq \liminf_{n\rightarrow\infty}\varphi_V(r_n) \leq  p\tilde{K}_V \tilde{A}_V < p\tilde{K}_V\tilde{B}_V \leq +\infty.$$
Thus we finish the proof.
\vskip2mm
\par
By using the similar arguments of Theorem 4.1, we can also obtain the similar results for the following scalar equation with Dirichlet boundary value on a locally finite graph $G=(V,E)$:
\begin{eqnarray}
\begin{cases}
\label{eqq2}
   \pounds_{m,p}u=\lambda f(x,u),\;\;\;\;\hfill x\in \Omega^\circ,\\
   |\nabla^j u|=0,\;\;\;\;\hfill x\in\partial\Omega, 0\leq j\leq m-1,
\end{cases}
\end{eqnarray}
where $p > 1, m\in \mathbb N$, $\lambda>0$  and $f:V\times\mathbb{R}\rightarrow\mathbb{R}$ and $\Omega\subset G(V,E)$ is a bounded domain.

\vskip2mm
  \noindent
{\bf Theorem 6.2.} {\it Suppose that $G=(V,E)$ is a locally finite graph, $F(x,s)=\int_0^s f(x,\tau)d\tau$ for all $x\in \Omega$, $\Omega^\circ\not=\emptyset $ and the following conditions hold: \\
 $(h)'$\: $h(x)>0$ for all $x\in \Omega$;\\
 $(f_0)'$\; $F(x,s)$ is continuously differentiable in $s\in \R$ for all $x\in \Omega$;\\
$(f_1)'$ \;   $\int_{\Omega}F(x,0)d\mu=0$;\\
$(f_2)'$ \;
$$0<\tilde {A}_\Omega:=\liminf_{y\rightarrow\infty}\frac{\int_{\Omega}\max_{|s|\leq y}F(x,s,t)d\mu}{y^p}<
 \limsup_{|s|\rightarrow\infty}\frac{\int_{\Omega} F(x,s)d\mu}{|s|^p}:=\tilde {B}_\Omega.$$
Then for each $\lambda\in(\tilde\lambda_{1,\Omega},\tilde\lambda_{2,\Omega})$ with $\tilde\lambda_{1,\Omega}=\frac{1}{\tilde{B}_\Omega}$ and $\tilde\lambda_{2,\Omega}=\frac{1}{p\tilde {K}_\Omega\tilde{A}_\Omega}$, where $\tilde {K}_\Omega=\frac{C^p(m,p,\Omega)}{\mu_{\min,\Omega}}$, equation (\ref{eqq2}) possesses an unbounded sequence of solutions.}

\vskip2mm
\par
By using the similar arguments of Theorem 5.1, we can also obtain the similar results for the following scalar equation on locally finite graph $G=(V,E)$:
\begin{eqnarray}
\label{eqq23}
  -\Delta_p u+h(x)|u|^{p-2}u=\lambda f(x,u),\;\;\;\;\hfill x\in V,
\end{eqnarray}
where $h:V\rightarrow\mathbb{R}$, $p\geq 2$, $\lambda>0$ and  $f:V\times\mathbb{R}\rightarrow\mathbb{R}$. We make the following assumptions:\\
$(h)$\; there exists a constant $h_0>0$ such that $h(x)\geq h_0>0$ for all $x\in V$;\\
$(M)$ there exists a $x_0\in V$ such that $M(x_0)\leq M(x)$ for all $x\in V$, where
$$
M(x)= \left(\frac{deg(x)}{2\mu(x)}\right)^{\frac{p}{2}}\mu(x)+h(x)\mu(x)+\sum_{y\thicksim x_0}\left(\frac{w_{xy}}{2\mu(y)}\right)^{\frac{p}{2}}\mu(y),\ \ x\in V.
$$
\par
Let
$$ \tilde \varrho= \frac{M(x_0)}{p} \;\;\mbox{and}\ \ \tilde K= \frac{1}{h_0^{1/p}\mu_0^{1/p}} .$$

\vskip2mm
  \noindent
{\bf Theorem 6.3.} {\it Let $G=(V,E)$ be a locally finite graph and $F(x,s)=\int_0^s f(x,\tau)d\tau$ for all $x\in V$. Assume that $(h)$, $(M)$ and the following conditions hold: \\
$(\tilde f_0)$\; $F(x,s)$ is continuously differentiable in $s\in \R$ for all $x\in V$,  and there exists a  function $a\in C(\R^+,\R^+)$ and a function $b:V\to \R^+$ with $b\in L^1(V)$ such that
$$
|F_s(x,s)|\le a(|s|) b(x),  |F(x,s)|\le a(|s|) b(x),
$$
for all $x\in V$ and all $s\in \R$;\\
 $(\tilde f_1)$\; $\int_V F(x,0)d\mu=0$;\\
 $(\tilde f_2)$ \;
$$0<\tilde A:=\liminf_{y\rightarrow\infty}\frac{\int_V \max_{|s|\leq y}F(x,s)d\mu}{y^p}<
\limsup_{|s|\rightarrow\infty}\frac{\int_V F(x,s)d\mu}{|s|^p}:=\tilde B.$$
Then for each $\lambda\in (\tilde\Theta_1,\tilde\Theta_2)$ with $\tilde\Theta_1=\frac{\tilde\varrho}{\tilde{B}}$ and $\tilde\Theta_2=\frac{1}{p\tilde {K}\tilde{A}}$, equation (\ref{eqq23}) possesses an unbounded sequence of solutions.}

 \vskip2mm
 \noindent
 {\bf Acknowledgments}\\
This project is supported by Yunnan Fundamental Research Projects of China  (grant No: 202301AT070465) and  Xingdian Talent Support Program for Young Talents of Yunnan Province of China.

\vskip2mm
 \noindent
 {\bf Authors' contributions}\\
The authors contribute the manuscript equally.

 \vskip2mm
 \noindent
 {\bf Competing interests}\\
The authors declare that they have no competing interests.

\vskip2mm
\renewcommand\refname{References}
{}
\end{document}